\documentclass[10pt]{amsart}

\usepackage[T1]{fontenc} 

\usepackage{amssymb,amsmath,amsfonts}

\usepackage{tikz}
\usepackage{algorithm}
\usepackage{algpseudocode}
\usepackage{float}
\usepackage{hyperref}
\usepackage{verbatim}

\usepackage{enumitem}

\newtheorem{theorem}{Theorem}
\newtheorem{lemma}{Lemma}

\newtheorem{definition}{Definition}

\newcommand{\lcm}{ {\rm lcm}}
\newcommand{\Z}{\mathbb{Z}}

\newtheorem{heuristic}[theorem]{Heuristic}

\makeatletter
\newtheorem{repeatthm@}{Theorem}

\makeatother

\algnewcommand\algorithmicforeach{\textbf{for each}}
\algdef{S}[FOR]{ForEach}[1]{\algorithmicforeach\ #1\ \algorithmicdo}

\title[Euler Pseudoprimes through Carmichael Products]{Constructing Large Euler Pseudoprimes through Carmichael Products} 

\author{Jon Grantham}
\author{Andrew Shallue}
\author{Jonathan Webster}

\address[1]{Institute for Defense Analyses/Center for Computing Sciences}
\address[2]{Illinois Wesleyan University}
\address[3]{Butler University}

\begin{document}

\begin{abstract}
Inspired by recent work on Carmichael numbers and by the construction in \cite{lange_etal}, we revisit the subset product construction of \cite{shallue2}.  
We state a theoretical parallel variant of the pushdown algorithm and 
extend computational records in two ways.
First, we construct a Carmichael number that is an Euler pseudoprime to every base up to $1172819$, which includes the first $90968$ prime numbers.
Second, we extend the known list of Carmichael numbers with exactly $d$ prime factors from $d \leq 19565220$ to $d \leq 340210163$.
\end{abstract}

\maketitle

\section{Introduction}

\textit{Carmichael numbers} are composite numbers $n$ for which $b^n \equiv b \pmod{n}$ for all $b \in \mathbb{Z}$.
They are pseudoprimes, since they satisfy the conclusion of Fermat's little theorem despite being composite.

These numbers have fascinated analytic number theorists and algorithmic number theorists alike, and the two groups have influenced each other in profound ways.
In the introduction to Alford, Granville, and Pomerance's proof of the infinitude of Carmichael numbers \cite{inf_carm}, their footnote credits the computational success of the subset-product problem as inspiring their own approach.
In the other direction, \cite{pom_self, dist_psp}'s arguments for the asymptotic count of Carmichael numbers provided the framework for Shallue \& Webster's proof of a heuristically optimal tabulation algorithm \cite[Theorem 3]{sw_ants4}.
We continue that tradition here, influenced by two recent results on Carmichael numbers.

\subsection{Carmichael numbers as Euler pseudoprimes}

Since Carmichael numbers are Fermat pseudoprimes, it has been common to see if they also satisfy the pseudoprimality of other, stronger tests based on Fermat's little theorem.  
One example is the Solovay-Strassen test, which given a candidate prime $n$ and a base $b$ with $\gcd(b,n)=1$ checks to see if 
\[ \left( \frac{b}{n} \right) \equiv b^{\frac{n-1}{2}} \pmod{n}  \]
and reports that $n$ is proven composite if the identity fails.  
The composite $n$ that pass the test are called Euler pseudoprimes to the base $b$.

In \cite{lange_etal} the authors present a novel construction.  
This involves finding a product of Carmichael numbers that forms an Euler pseudoprime to all integer bases $b$ up to $211$, i.e. it fools the Solovay-Strassen test for the smallest $47$ prime bases.   
Boltnev \cite{boltnev} previously constructed a strong (and thus Euler) pseudoprime to the first $350$ prime bases, through $2357$. 
His example also passes the next two prime bases, so it covers all prime bases up to $2377$.
We present a new construction which greatly improves upon this result; we have found a product of Carmichael numbers that forms an Euler pseudoprime to all integer bases $b$ up to $1172819$. 

\subsection{Carmichael numbers with  \texorpdfstring{$d$}{}  prime factors}
It had long been conjectured that for every $d > 2$, there exists a Carmichael number with exactly $d$ prime factors.  
Larsen and Wright recently proved an asymptotic version of this result showing that for all $d$ sufficiently large, there exists a Carmichael number with exactly $d$ prime factors \cite{larsen-wright}.  
In private correspondence, Larsen indicated that it might be possible to make the bound effective and guessed that doing so would result in $d$ on the order of $10^{1000}$.  
This is significantly beyond our ability to bridge the known counts of $19565220$  from \cite{shallue2}.
Nonetheless, we revisited this computation and extended this result to $340210163$, a factor of $17$ improvement.

\subsection{Organization}

The remainder of the paper is structured as follows.  
Section \ref{sec:background} provides the necessary background on Carmichael numbers and pseudoprimes and provides a summary of the pseudoprime construction technique of \cite{lange_etal}.
Section \ref{sec:subset_prod} reviews Algorithm 1 of \cite{shallue2} and explains how it may be run efficiently in a parallel computing environment. 
Section \ref{sec:epsp} explains the required modifications and additions to the subset product algorithm to find an Euler pseudoprime to many bases.  
Section \ref{sec:implementation} provides some details of the implementation and of the computational results.

\section{Background} \label{sec:background}

We frequently recognize Carmichael numbers $n$ via their prime factorization.

\begin{theorem}[Korselt's Criterion]
A composite number $n$ is a Carmichael number if and only if $n$ is square-free and $(p-1) \mid (n-1)$ for all prime divisors $p$ of $n$. 
\end{theorem}

Robert Carmichael proved Korselt's criterion in the context of the reduced totient function or the Carmichael function $\lambda(n)$, 
\begin{align*}
\lambda\left(p^{\alpha}\right) &= \begin{cases} \phi\left(p^{\alpha}\right) &\text{ if } \alpha \leq 2 \text{ or } p \geq 3 \\ \frac{1}{2} \phi\left(p^{\alpha}\right) &\text{ if } \alpha \geq 3 \text{ and } p = 2 \end{cases},
\end{align*}
\[ \lambda(n) = \lambda\left( {p_1}^{\alpha_1} \cdots {p_k}^{\alpha_k} \right) = \lcm \, \big(\lambda\left({p_1}^{\alpha_1}\right), \ldots, \lambda\left({p_k}^{\alpha_k}\right)\big), \]
where the $p_i$ are distinct prime divisors of $n$ and $\phi(n)$ is Euler's totient function.  
Later we will also refer to $v_2(n)$, the $2$-adic valuation of $n$, which returns the largest exponent $e$ such that $2^e \mid n$.
We can now state Korselt's criterion in terms of $\lambda(n)$.

\begin{theorem}[Korselt's Criterion]
A composite number $n$ is a Carmichael number if and only if $\lambda(n) \mid (n-1)$. 
\end{theorem}

Some early views on the infinitude of Carmichael numbers came from the fact that $n = (6k + 1)(12k + 1)(18k + 1)$ is a Carmichael number if $n$ is a product of three primes.  
If the prime $k$-tuples conjecture were true, then this family would produce infinitely many Carmichael numbers.
However, the argument that won the day (often attributed to Erd\H{o}s) is centered not around specific families of $n$, but rather around numbers with some prescribed Carmichael Lambda value.  
Indeed, let
\[ \Lambda = \prod p_{i}^{e_i},\]
where $p_i$ are the ordered prime numbers and $e_i > 0$.  
Consider the set $\mathcal{P} = \{ p \mbox{ prime: }  (p - 1) \mid \Lambda \ , \  p \nmid \Lambda \}$.
If $2^{|\mathcal{P}|}$ is large compared to $\phi( \Lambda )$, then it should be possible to find a subset product in $\mathcal{P}$ that is $1 \pmod{\Lambda}$. 
The idea is that a subset product in $\mathcal{P}$ takes a random value in $(\mathbb{Z}/\Lambda\mathbb{Z})^\times$ and with enough such products it is likely that at least one such composite product will land on $1$.
Let $N$ be such a product.
Then $\lambda(N) \mid \Lambda \mid (N - 1)$, and so it is a Carmichael number.
It is possible to choose $\Lambda$ so that many pair-wise co-prime Carmichael numbers may be found as subset products in $\mathcal{P}$, and this was the key strategy of \cite{shallue2}, which we explain in Section \ref{sec:subset_prod}.

\subsection{Constructing an Euler pseudoprime from Carmichael numbers}

See \cite{lange_etal} for a fuller treatment of this topic.

\begin{definition}[Index]
The {\it index} of a Carmichael number $n$ is given by 
\[ i(n) = \frac{n-1}{\lambda(n)} \enspace . \]
\end{definition}

In \cite{lange_etal}, a coarse classification of Carmichael numbers was constructed. 
We focus only on the class needed for this paper.

\begin{definition}[Class A Carmichael numbers]
A Carmichael number $n$ is considered to be in Class A if $i(n)$ is even.
\end{definition}

\begin{lemma}[\cite{lange_etal}, Lemma 7] \label{lemma1}
Let $n$ be a Carmichael number in Class A.  Then there are $\frac{1}{2} \phi(n)$ Euler liars in $\Z_n^\times$ and all have Jacobi symbol $1$.
\end{lemma}

The Class A Carmichael numbers make good candidates for Euler pseudoprimes.
We need a Lemma regarding the product of two Class A Carmichael numbers.

\begin{lemma}[\cite{lange_etal}, Lemma 12] \label{lemma2}
Let $n_1$ and $n_2$ with $\gcd (n_1, n_2) = 1$ be Class A Carmichael numbers.
If $n = n_1n_2$ is a Carmichael number, then $n$ is in Class A if $\min \{ v_2( n_1 -1), v_2(n_2 -1) \} > \max \{ v_2( \lambda(n_1)), v_2( \lambda(n_2) ) \}$.
\end{lemma}

So, these two lemmas allow us to multiply Class A Carmichael numbers and retain their classification and Euler pseudoprimality properties.  

\begin{lemma}[\cite{lange_etal}, Lemma 13]\label{lemma3}
Let $n_1$ and $n_2$ be Class A Carmichael numbers that are Euler pseudoprimes for every base less than the prime $b$ and not to $b$ itself.  
Assume $\gcd(n_1, b) = \gcd(n_2, b) = 1$. 
Then if $n = n_1n_2$ is a Class A Carmichael number, it is an Euler pseudoprime for all bases up to and including $b$.
\end{lemma}

\begin{proof}
Assume $n_1$ and $n_2$ are in  Class A, which means that the multiplicative order of $a$ modulo $n_i$ for all $a < b$ is a divisor of $(n_i-1)/2$.
Thus the Jacobi symbol is $1$, and so the Jacobi symbol for the product $n_1 n_2$ is also $1$ for all bases less than $b$.  
This makes the product also an Euler pseudoprime to all bases less than $b$.

Now focusing on $b$, the Jacobi symbol with $b$ is $-1$.  
By multiplying two such numbers, the Jacobi symbol becomes $1$, and the product being Class A means that the power is also $1$.  
This makes the product an Euler pseudoprime to the base $b$ as well.
\end{proof}

The above is sufficient to describe their method and we refer the readers to Sections 3--5 of \cite{lange_etal} for their complete description.
They used the table of Carmichael numbers produced by \cite{sw_ants4} to construct products that were again Carmichael numbers.
They repeatedly invoked Lemma \ref{lemma3} in order to construct Euler pseudoprimes to consecutive prime bases.
At each stage of their algorithm, they were able to account for one more prime base.

The practical difficulty in their construction is the antecedent ``if $n = n_1n_2$ is a Carmichael number.'' 
This is unlikely to hold.
Our solution, alluded to above, is to use the subset product construction to guarantee the antecedent is always satisfied.  

\begin{lemma}\label{lem:product_is_car}
Let $\Lambda$ be a positive integer.  Assume that $n_1$ and $n_2$ are Carmichael numbers with $\gcd(n_1, n_2) = 1$.  Further assume
\[
(p-1) \mid \Lambda \mid (n_1 - 1) \hspace{2mm} \forall \ p \mid n_1   
\mbox{   and   }
(p-1) \mid \Lambda \mid (n_2 - 1) \hspace{2mm} \forall \ p \mid n_2.
\]
Then $n_1 \cdot n_2$ is also Carmichael.
\end{lemma}
\begin{proof}
Since $\gcd(n_1, n_2) = 1$ and both are Carmichael numbers, their product is squarefree.  Since $n_1 \equiv 1 \pmod{\Lambda}$ and $n_2 \equiv 1 \pmod{\Lambda}$, their product is also equivalent to $1$ modulo $\Lambda$.
The divisibility assumption of the Korselt criterion is then satisfied by the divisibility assumption of the lemma.
\end{proof}

In the next section we will see that the result of Algorithm \ref{alg:subsetprod} naturally leads to the conditions of Lemma \ref{lem:product_is_car} being satisfied, and hence our Carmichael products will always themselves be Carmichael.

The only minor adjustment, to be explained below, is to guarantee the inequality on the 2-valuations in Lemma \ref{lemma2} is always satisfied.
Having done that, we will show that all of the Jacobi symbols may be aligned at once rather than one at a time.  
For our purposes, a simpler version of Lemma \ref{lemma3} suffices, which we state without proof.  

\begin{lemma}\label{lemma4}
Let $n_1$ and $n_2$ be Class A Carmichael numbers that are not base $b$ Euler pseudoprimes. 
Assume $\gcd(n_1, b) = \gcd(n_2, b) = 1$. 
Then if $n = n_1n_2$ is a Class A Carmichael number, it is a base $b$ Euler pseudoprime.
\end{lemma}

\section{The Subset Product Algorithm}\label{sec:subset_prod}

We present a slight modification to the modular subset product algorithm of \cite{shallue2} suitable for a parallel computing environment.  
As set up, we choose some 
\[ \Lambda = \prod_{i=1}^r p_{i}^{e_i},\]
where $p_i$ are the ordered prime numbers and $e_i > 0$.  
We let $Q_i = p_i^{e_i}$; in particular $Q_1$ is some power of $2$.
We construct the set $\mathcal{P} = \{ p \mbox{ prime: }  (p - 1) \mid \Lambda \ , \  p \nmid \Lambda \}$, then find many disjoint subsets of $\mathcal{P}$ that product to $1$ modulo $\Lambda$.  
Each such composite product is a Carmichael number.

\subsection{Asymptotic Complexity}

There is a rich literature on the modular subset product problem, only some of which we will review here.  
We start with the observation that the problem is NP-complete \cite{karp_npcomplete}.  
Despite this fact, there are many instances where a polynomial time algorithm exists.  
Understanding which instances are easy and which are hard requires the notion of density.

\begin{definition}
Let $\{x_1, x_2, \dots, x_n, b, M\}$ define a modular subset product problem, 
where our goal is to find a subset of the $x_i$ that products to $b$ modulo $M$.

The density of the problem is defined as $\frac{n}{\log{\phi(M)}}$, where $\log$ is the base-$2$ logarithm.
\end{definition}

In a rough sense, the density compares the count of subsets with the size of the group.  
For a random instance of the problem, a solution exists with high probability if and only if the density is greater than one.  
If the density is less than one and we are confident that a solution exists, a good strategy is to use lattice basis reduction \cite{lowdensity_subsetsum}.  
If the density is one (the hardest instances), one could use a birthday-problem approach that costs $\widetilde{O}(2^{n/2})$ time and space, or an improved version like \cite{beckercoronjoux} with running time $\widetilde{O}(2^{0.291n})$ for additive subset sum.

When constructing Carmichael numbers it is easy to generate instances with density much greater than one, so we will now focus on that case.  
The classic solution is dynamic programming: build a 2D table and fill it in at a cost of $O(n\cdot M)$ time and space.  
In our case $M$ is large enough that this solution becomes impractical, so we turn instead to the method of \cite{lohniebuhr}, which was improved in \cite{shallue2}.  
See the next subsection for the presentation and discussion of that algorithm.

\subsection{Algorithm}

In our construction of an Euler pseudoprime to many bases we needed to construct many Carmichael numbers, all with respect to the same $\Lambda$ and the same set $\mathcal{P}$.  
For this purpose we used \cite[Algorithm 1]{shallue2}.  

We call this the ``pushdown algorithm,'' because at each step pairs of products are found which are multiplicative inverses modulo the next $Q_u$, pushing the new product down to a smaller subgroup of $(\Z/\Lambda\Z)^{\times}$.
By repeating the pushdown step for each of the prime powers $Q_u$ dividing $\Lambda$, the algorithm ends with a collection of disjoint products, all of which are the identity modulo $\Lambda$, and hence each composite product is a Carmichael number.

The algorithm requires high density in order to succeed.  Furthermore, choosing a $\Lambda$ with higher density results in more Carmichael numbers constructed.

\begin{algorithm}[htbp]
\caption{Many Carmichael numbers}
\label{alg:subsetprod}
\begin{algorithmic}[1]
\State $S_0 \leftarrow \mathcal{P}$

\For{$u \leftarrow 1$ \textbf{to} $r$} 
    \State $S_u \leftarrow \emptyset$
    \State (Optional) $\mathrm{fixable} \leftarrow \mathrm{true}$
    \State (Optional) Calculate $b_u = \prod_{a \in S_{u-1}} a \pmod{Q_u}$
    \If{(Optional) $b_u \ne 1$}
        \If{there exists $a \in S_{u-1}$ with
            $a \equiv b_u \pmod{Q_u}$}
            \State Remove exactly one such $a$ from $S_{u-1}$
        \Else
            \State $\mathrm{fixable} \leftarrow \mathrm{false}$
        \EndIf
    \EndIf
    \For{$a \in S_{u-1}$}   
        \If{$a\equiv 1 \pmod{Q_u}$ }
            \State $S_{u} \leftarrow S_u \cup \{a\}$ and remove $a$ from $S_{u-1}$
        \Else
            \State Find $\hat a \in S_{u-1} \setminus \{a\}$ such that $a\hat a \equiv 1 \pmod{Q_u}$, if one exists
            \If{such an $\hat a$ exists}
                \State $S_u \leftarrow S_u \cup \{a\hat a\}$, and remove $a, \hat a$ from $S_{u-1}$
            \EndIf
        \EndIf
    \EndFor 
    \If{(Optional) $\mathrm{fixable}$ and $S_{u-1} \ne \emptyset$}
        \State (Optional)  Multiply the remaining $a\in S_{u-1}$, add to $S_{u}$
    \EndIf
\EndFor
\State \Return  $S_{r}$, each composite element of which is Carmichael  
\end{algorithmic}
\label{alg:many_cars}
\end{algorithm}

Lines 4--12 and 23--25 are optionally included (or excluded) all at once.  
Not executing the optional lines decreases the number of elements that would have been in $S_u$ by at most $1$ per round for a given $S_{u-1}$. 
Since the goal is to have many Carmichael numbers,  we assume the optional lines are executed, and this is how the original Algorithm 1 of \cite{shallue2} worked.
We address the significance of these optional lines in Subsection \ref{subsec:parallel}.

We note that lines 13--22 can be accomplished with a bucket-style sort using at most $\phi(Q_u)/2+1$ buckets. 
Each bucket is jointly used by $a$ or $\hat a$ and can only have one or the other element (if both are found, it is a match and removed).  

An important implementation detail is that group elements $a$ are stored as a list of integers which are elements of $\mathcal{P}$ (called the history), along with their product as an element of $(\Z/\Lambda\Z)^{\times}$.  
In this way we can reconstruct the Carmichael number as an integer, while being able to calculate the needed residues in a constant number of group operations.

\begin{theorem}\label{thm:determ_complexity}
Algorithm \ref{alg:subsetprod} has time complexity $O(rn)$, measured in operations in the group $(\Z/\Lambda\Z)^{\times}$.  Here $r$ is the number of distinct primes dividing $\Lambda$ and $n = |\mathcal{P}|$.
\end{theorem}
\begin{proof}
The algorithm has a nested loop structure, where the outer loop has length $r$ and the inner loop has length $|S_{u-1}| \leq |S_0| = n$.

For a given $u$, processing via the bucket-style sort requires one pass at a cost of $O(|S_{u-1}|)$ group operations, since for each $a$ we either pass it along, put it in a bucket, or pair it with an element in a bucket and pass the pair along.
\end{proof}

\begin{heuristic}\label{thm:expected_complexity}
Assume the elements of $\mathcal{P}$ are uniformly distributed elements of $(\Z/\Lambda\Z)^{\times}$.  Then
Algorithm \ref{alg:subsetprod} has expected time complexity $O(|\mathcal{P}|)$, measured in group operations.
\end{heuristic}
\begin{proof}[Heuristic justification]
The total number of $a$ considered by the algorithm is 
$$
|S_0| + |S_1| + \cdots + |S_{r-1}| \enspace .
$$
At each step in the outer loop, an expected fraction $\frac{1}{\phi(Q_u)}$ of the elements will be passed to the next stage immediately.  Of the remaining, we expect them to pair up and for their product to be passed along.  Thus we have the expectation that 
$$
|S_{u}| \leq \frac{|S_{u-1}|}{\phi(Q_u)} + \frac{1}{2}\left(|S_{u-1}| - \frac{|S_{u-1}|}{\phi(Q_u)} \right)
\leq |S_{u-1}| \left( \frac{1}{2 \phi(Q_u)} + \frac{1}{2} \right)
\leq \frac{2}{3} |S_{u-1}|
$$
where the last inequality follows as long as $\phi(Q_u) \geq 3$. (At most two rounds fail that inequality, and they do not affect the overall bound.)
We conclude that the total number of $a$ is $O(|S_0|) = O(|\mathcal{P}|)$ by bounding the resulting geometric sum.

At step $u$ we scan through the set of elements, either placing an element in one of the $\phi(Q_u)/2+1$ buckets or matching with an element in a bucket and removing both.  
Since we have access to the residue modulo $\Lambda$, processing each element results in a constant number of group operations.
\end{proof}

Note that counting bit operations rather than group operations would make Theorem \ref{thm:determ_complexity} a more accurate picture compared to Heuristic \ref{thm:expected_complexity}, as effectively the entirety of $S_0$ is carried along from step to step via the history portion of the element representation.

These bounds make no claims as to the density required to ensure that $S_{r}$ is nonempty with high probability; in practice we found it easy to construct $\Lambda$ with high enough density.

\subsection{A parallel version of Algorithm~\ref{alg:many_cars}} \label{subsec:parallel}

We describe two versions of a theoretical distributed-memory implementation of
Algorithm~\ref{alg:many_cars}.
The key distinction is found in the inclusion or exclusion of the optional lines.
The first, where we exclude the optional lines, is more elementary and provides a simple chunking structure.  
This version may actually be run sequentially if it is not feasible to have all of $\mathcal{P}$ stored at once.
It also uses no inter-process communication and may be preferred if inter-process communication is especially costly. 
The second version shows how to introduce inter-process communication and allow the execution of the optional lines.

\subsubsection*{Model}
We assume $k$ processes with ranks $0,\dots,k-1$, each with private memory, and
point-to-point and tree-structured collective communication.

\subsubsection*{Parallel prime generation}
The ability to generate primes in parallel was addressed in \cite{shallue2}.
One simple method to accomplish this is by assigning certain fixed divisors of $\Lambda$ to given processors.  
In this way, we can create
\[ \mathcal{P} = \bigcup_{j = 0}^{k-1} S_0^{(j)}, \]
where the union is a disjoint union and processor $j$ has $S_0^{(j)}$.
With the exception of the assignment of the fixed divisors for initialization, prime generation needs no communication among the $k$ processors.  

\subsubsection*{Trivial parallelization}
Processor $j$ simply runs Algorithm \ref{alg:many_cars}, where line $1$ has the assignment $S_0 \leftarrow S_0^{(j)}$.  
After each of $k$ processors finish their work, the elements in $S_{r}^{(j)}$ can be pooled together.
The success of Algorithm \ref{alg:many_cars} relies on having a sufficiently dense set of solutions to the subset product problem. 
This places a practical limit on the number of processors that can be used in this trivial parallelization before the expectation is that each local job would fail to find Carmichael 
numbers.  
Conservatively, each process might need
\[ 
|S_0^{(j)}| \ \gtrsim\ \max_{1\le u\le r} 2^{u-1} Q_u^2 ,
\]
primes so that the final step is expected to be nonempty.  
We can use static partitioning by choosing 
\[ k \ \approx\ \frac{|\mathcal{P}|}{ \max_{1\le u\le r} 2^{u-1}Q_u^2} . \]
Heuristic estimates of $|\mathcal{P}|$ are good once $\Lambda$ is fixed, so $k$ can be
chosen before the computation begins.

\subsubsection*{Parallel version with communication}

We explain how to accomplish Algorithm \ref{alg:many_cars} in parallel while executing the optional lines.  
The guiding principle is that the elements of $S_u$ grow as $u$ increases, but the information needed to decide which elements to combine is small (residues modulo $Q_u$).
The algorithm therefore favors local work, communicates residues modulo $Q_u$ freely, and communicates the elements $a$ or $\hat a$ as rarely as possible.
It does not matter if $a$ is treated as an actual integer formed as a product of primes or as a data structure that stores the prime factors explicitly.  
In either case, the asymptotic cost of communication remains the same.

We describe round $u$ of Algorithm \ref{alg:many_cars} in a parallel environment.
As with the trivial parallelization, line $1$ is replaced with $S_0 \leftarrow S_0^{(j)}$.  
And rank $j$ only has access to $S_{u-1}^{(j)}$ unless the instructions indicate communication.
\begin{enumerate}
    \item Line 5 is computed by an all-reduce of the residues.
    \item For lines 6--12, when $b_u \ne 1$, each process reports whether it holds an element congruent to $b_u \pmod{Q_u}$ and the lowest such rank, found with an all-reduce, deletes one such element. If none exist, set fixable to false.
    \item Lines 13--22 are first accomplished locally.  Unmatched elements are saved for parallel matching. 
    \item The remaining matches are found by a binary tree reduction over the residue tables.  
    \begin{itemize}
        \item Pairs of processors share their residue tables which contain the rank and surplus residues.  When a match is found, only $a$ or $\hat a$ are communicated as necessary. 
        \item The pairs produce a joint residue and rank table, which is used in the node above.
    \end{itemize}
    \item For lines 23--25, if fixable is true and the global leftover pool is nonempty, once the root's table has no more matches, each process multiplies its own unmatched elements. These $k$ partial products are combined in a tree product into one last element of $S_{u}$, which may be assigned to any processor.
\end{enumerate}


The intent is to write this parallel code using a static partitioning of the work.  
Keeping in mind that most work is meant to be done locally, it is possible to redesign the above to allow dynamic partitioning.
If $S_u^{(j)}$ drops below some threshold, processor $j$ can send all of its products in the first step of the tree reduction.    

\section{Euler Pseudoprimes to many bases --- Linear Algebra over \texorpdfstring{$\mathbb{F}_2$}{} } \label{sec:epsp}

It is simple to modify the above to create Class A Carmichael numbers whose products are also always in Class A.
We need the inequality on the 2-valuations of Lemma \ref{lemma2} to be satisfied.  
There are two conceptually equivalent ways to accomplish this.  
First, we can run the prime generation phase on $\Lambda/2$.  
Carmichael numbers produced this way satisfy $\lambda(n) \mid (\Lambda/2) \mid \Lambda \mid (n-1)$.
Second, we can replace $Q_1$ (the power of $2$) with $2Q_1$ when we run Algorithm \ref{alg:many_cars}.
Numbers produced in this way satisfy $\lambda(n) \mid \Lambda \mid 2\Lambda \mid (n-1)$.
We either artificially lower the power of $2$ for the prime generation phase or we artificially inflate the power of $2$ for the subset product problem.  
In either case, $v_2(n-1) > v_2(\lambda(n) )$, and for the remainder of the section we assume one of these two approaches has been taken.  

Consider a set $\mathcal{B}$ of bases and let $b_k$ denote the $k$-th element of $\mathcal{B}$.  The easiest case to consider (and the case we computed) involves consecutive primes starting with $2$, but in principle $\mathcal{B}$ could be any set of positive integers.  For Lemma \ref{lemma4} to apply we need $\gcd(n, b) = 1$ for all $n \in S_r$ and all $b \in \mathcal{B}$.  In our computation we achieved this by restricting the primes in $\mathcal{P}$ before running the pushdown algorithm to ensure all were greater than the largest element of $\mathcal{B}$.

By Lemma \ref{lemma1}, we know half of all coprime bases are Euler liars (and thus have symbol $1$).  
For each $n_i \in S_{r}$ associate to it a vector $ \textbf{b}_i$ over $\mathbb{F}_2$ in the following way.
Let the coordinate in the $k$-th place of $\textbf{b}_i$ have the value $0$ if $b_k$ is an Euler liar for $n_i$ and $1$ otherwise (that is, the Jacobi symbol would be $-1$).
So long as $|S_{r}| > |\mathcal{B}|$, we can use linear algebra over $\mathbb{F}_2$ to find a nontrivial solution $\textbf{x} $ to
\[ [ \textbf{b}_1 \ \textbf{b}_2 \ \cdots \ \textbf{b}_{|S_{r}|} ] \ \textbf{x} = \textbf{0}.\]
The solution corresponds to a disjoint product of Carmichael numbers in $S_{r}$ which is an Euler pseudoprime to every $b_k \in \mathcal{B}$.  This product is itself Carmichael by Lemma \ref{lem:product_is_car}.

\section{Implementation Details}\label{sec:implementation}

The code for this paper and the factors of the Euler pseudoprime are available at \url{https://github.com/31and8191/euler-carmichael-computations}.

\subsection{Euler pseudoprime computation}

Let 
\begin{align*}
\Lambda &=2^8 \cdot 3^5 \cdot5^3 \cdot7^2 \cdot11^2 \cdot13^2
\prod_{17\le p \le73}p\\
       &=10567690344536672406900888090467424000.
\end{align*}
For $d \mid \Lambda$, we look for primes $p=d+1$. We require $p$ to be greater than $1.3$ million, so that we do not have the Jacobi symbol equal to zero for any tested base. By using $2Q_1=2^9$ for pushdown, we ensure the resulting Carmichael numbers are Class A. The computation found $19551404$ primes which combined via pushdown to give $90943$ Carmichael numbers. The prime generation took less than 4 minutes with 24 jobs, and the serial pushdown took less than 9 minutes.

We form the matrix described in the previous section whose columns correspond to Carmichael numbers and rows to prime bases not dividing $\Lambda$. The matrix with $90943$ columns and $90947$ rows has a one-dimensional kernel. We used a Python implementation of Gaussian elimination over $\mathbb{F}_2$ to find a nonzero solution. Along with the $21$ bases dividing $\Lambda$, the resulting Carmichael number passes the Euler pseudoprime test for the first $90968$ prime bases.

The solution vector corresponds to the product of $45413$ of the Carmichael numbers, with a total of $8667018$ prime factors. This number has $160558499$ decimal digits. Constructing the final Carmichael number from the $90943$ seed Carmichael numbers, including the matrix work, took approximately 11.5 hours.

We then used a separate program to verify that the factors are distinct, for each $p$, $(p-1) \mid \Lambda$, and that the product $n$ has $2\Lambda \mid (n-1)$. For each prime base up to $1172819$, we computed the Legendre symbols $(a/p)$ and multiplied them to get the Jacobi symbol $(a/n)$. It was $1$ for every tested base. Because $n$ is Class A, we get the Euler condition.

\subsection{Carmichael numbers with \texorpdfstring{$d$}{} prime factors}

With a new implementation of Algorithm 1 of \cite{shallue2}, we extended the computation in that paper to demonstrate the existence of Carmichael numbers with exactly $d $ prime factors for every $3 \le  d \le 340210163$.

We use
\[
  \Lambda=2^{12}\cdot  3^7 \cdot 5^4 \cdot 7^3 \cdot 11^2 \cdot  13^2 \cdot  17^2 
  \prod_{\substack{19\leq p\leq83}}p.
  \]
We obtain $340210290$ primes which we combine via pushdown into $1529840$ base Carmichael numbers. We compute the possible numbers of prime factors of products of these Carmichael numbers and see that there are Carmichael numbers with $d$ factors for all $32\le d \le 340210163$. Combined with the original result, this gives all  $3\le d \le 340210163$.

The prime generation took about 1 hour and 42 minutes with 24 jobs, the pushdown took about 37 hours, and checking the base Carmichael numbers took about 7 minutes.

\section{Conclusion}

We are indebted to the authors of \cite{lange_etal} for an inspiring paper.  
Among other reasons, this gave us a perfect excuse to revisit the pushdown algorithm, eventually parallelize it, and create a new codebase for it.

Note that we took the prime bases up to some bound for convenience, but the techniques can be applied to any chosen set of bases.  It seems likely that similar techniques could be applied to the construction of other pseudoprimes.

\section*{AI Disclosure}

The authors used Anthropic Claude Opus 4.8 to write code for this project. They used OpenAI Codex for subsequent code development and review, checking numerical details, and paper planning. The authors wrote the text of the paper, with  editorial support provided by LLMs.  The authors have checked the content of this paper and take responsibility for it.

\bibliographystyle{amsplain}
\bibliography{manybases}

@article {inf_carm,
    AUTHOR = {Alford, W. R. and Granville, Andrew and Pomerance, Carl},
     TITLE = {There are infinitely many {C}armichael numbers},
   JOURNAL = {Ann. of Math. (2)},
  FJOURNAL = {Annals of Mathematics. Second Series},
    VOLUME = {139},
      YEAR = {1994},
    NUMBER = {3},
     PAGES = {703--722},
      ISSN = {0003-486X},
     CODEN = {ANMAAH},
   MRCLASS = {11N13 (11A07 11A51 11N25)},
  MRNUMBER = {1283874 (95k:11114)},
MRREVIEWER = {H. L. Montgomery},
       DOI = {10.2307/2118576},
       URL = {http://dx.doi.org/10.2307/2118576},
}

@article {pom_self,
    AUTHOR = {Pomerance, Carl and Selfridge, J. L. and Wagstaff, Jr., Samuel
              S.},
     TITLE = {The pseudoprimes to {$25\cdot 10^{9}$}},
   JOURNAL = {Math. Comp.},
  FJOURNAL = {Mathematics of Computation},
    VOLUME = {35},
      YEAR = {1980},
    NUMBER = {151},
     PAGES = {1003--1026},
      ISSN = {0025-5718},
     CODEN = {MCMPAF},
   MRCLASS = {10A40 (10-04 10A25)},
  MRNUMBER = {572872 (82g:10030)},
       DOI = {10.2307/2006210},
       URL = {http://dx.doi.org/10.2307/2006210},
}

@article {shallue2,
    AUTHOR = {Alford, W. R. and Grantham, Jon and Hayman, Steven and
              Shallue, Andrew},
     TITLE = {Constructing {C}armichael numbers through improved
              subset-product algorithms},
   JOURNAL = {Math. Comp.},
  FJOURNAL = {Mathematics of Computation},
    VOLUME = {83},
      YEAR = {2014},
    NUMBER = {286},
     PAGES = {899--915},
      ISSN = {0025-5718},
   MRCLASS = {11Y16 (11Y05)},
  MRNUMBER = {3143697},
MRREVIEWER = {Emil Vatai},
       DOI = {10.1090/S0025-5718-2013-02737-8},
       URL = {http://dx.doi.org/10.1090/S0025-5718-2013-02737-8},
}

@article {dist_psp,
    AUTHOR = {Pomerance, Carl},
     TITLE = {On the distribution of pseudoprimes},
   JOURNAL = {Math. Comp.},
  FJOURNAL = {Mathematics of Computation},
    VOLUME = {37},
      YEAR = {1981},
    NUMBER = {156},
     PAGES = {587--593},
      ISSN = {0025-5718,1088-6842},
   MRCLASS = {10A21 (10A20)},
  MRNUMBER = {628717},
MRREVIEWER = {M.\ C.\ Wunderlich},
       DOI = {10.2307/2007448},
       URL = {https://doi.org/10.2307/2007448},
}

@misc{lange_etal,
  author = {Alcantarilla S{\'a}nchez, Alejandra and
            Cottaar, Jolijn and
            Lange, Tanja and
            de Weger, Benne},
  title = {Ours go to 211: {Euler} pseudoprimes to 47 prime bases
           (from {Carmichael} numbers)},
  year = {2026},
  note = {\href{https://arxiv.org/abs/2602.21840}{arXiv:2602.21840}}
}

@article{boltnev,
 author = {Boltnev, Yu. F.},
 title = {Construction of strong pseudoprime {Carmichael} numbers passing the {Miller}-{Rabin} test for many first prime bases},
 fjournal = {Prikladnaya Diskretnaya Matematika. Supplement},
 journal = {Prikl. Diskretn. Mat., Suppl.},
 issn = {2226-308X},
 number = {18},
 pages = {258--262},
 year = {2025},
 language = {Russian},
 zbMATH = {8120180},
 Zbl = {1573.11015}
}

@article{larsen-wright,
 author = {Larsen, Daniel and Wright, Thomas},
 title = {Carmichael numbers with a specified number of prime factors},
 fjournal = {Integers},
 journal = {Integers},
 issn = {1553-1732},
 volume = {26A},
 pages = {A15},
 note = {14 pp.}, year = {2026},
 doi = {10.5281/zenodo.23018145},
 zbMATH = {8268338}
}

@inproceedings{sw_ants4,
  author = {Shallue, Andrew and Webster, Jonathan},
  title = {Algorithms for {C}armichael numbers},
  booktitle = {Proceedings of the Seventeenth Algorithmic Number Theory Symposium (ANTS XVII)},
  year = {2026},
  note = {To appear; \url{https://arxiv.org/abs/2506.09903v3}},
}

@incollection {karp_npcomplete,
    AUTHOR = {Karp, Richard M.},
     TITLE = {Reducibility among combinatorial problems},
 BOOKTITLE = {Complexity of computer computations ({P}roc. {S}ympos., {IBM}
              {T}homas {J}. {W}atson {R}es. {C}enter, {Y}orktown {H}eights,
              {N}.{Y}., 1972)},
    SERIES = {The IBM Research Symposia Series},
     PAGES = {85--103},
 PUBLISHER = {Plenum, New York-London},
      YEAR = {1972},
   MRCLASS = {68A20},
  MRNUMBER = {378476},
MRREVIEWER = {John\ T.\ Gill},
}

@article {lowdensity_subsetsum,
    AUTHOR = {Coster, Matthijs J. and Joux, Antoine and LaMacchia, Brian A.
              and Odlyzko, Andrew M. and Schnorr, Claus-Peter and Stern,
              Jacques},
     TITLE = {Improved low-density subset sum algorithms},
   JOURNAL = {Comput. Complexity},
  FJOURNAL = {Computational Complexity},
    VOLUME = {2},
      YEAR = {1992},
    NUMBER = {2},
     PAGES = {111--128},
      ISSN = {1016-3328,1420-8954},
   MRCLASS = {11Y16 (68Q25)},
  MRNUMBER = {1190825},
       DOI = {10.1007/BF01201999},
       URL = {https://doi.org/10.1007/BF01201999},
}

@article {lohniebuhr,
    AUTHOR = {L\"oh, G\"unter and Niebuhr, Wolfgang},
     TITLE = {A new algorithm for constructing large {C}armichael numbers},
   JOURNAL = {Math. Comp.},
  FJOURNAL = {Mathematics of Computation},
    VOLUME = {65},
      YEAR = {1996},
    NUMBER = {214},
     PAGES = {823--836},
      ISSN = {0025-5718,1088-6842},
   MRCLASS = {11Y11 (11A51 11Y16)},
  MRNUMBER = {1325872},
MRREVIEWER = {Andrew\ Granville},
       DOI = {10.1090/S0025-5718-96-00692-8},
       URL = {https://doi.org/10.1090/S0025-5718-96-00692-8},
}

@incollection {beckercoronjoux,
    AUTHOR = {Becker, Anja and Coron, Jean-S\'ebastien and Joux, Antoine},
     TITLE = {Improved generic algorithms for hard knapsacks},
 BOOKTITLE = {Advances in cryptology---{EUROCRYPT} 2011},
    SERIES = {Lecture Notes in Comput. Sci.},
    VOLUME = {6632},
     PAGES = {364--385},
 PUBLISHER = {Springer, Heidelberg},
      YEAR = {2011},
      ISBN = {978-3-642-20464-7},
   MRCLASS = {94A60 (68Q25)},
  MRNUMBER = {2813651},
       DOI = {10.1007/978-3-642-20465-4\_21},
       URL = {https://doi.org/10.1007/978-3-642-20465-4_21},
}

\end{document}